\documentclass{amsproc}%
\usepackage{amsfonts}
\usepackage{amsmath}
\usepackage{amssymb}
\usepackage{graphicx}
\usepackage{hyperref}%
\providecommand{\U}[1]{\protect\rule{.1in}{.1in}}
\theoremstyle{plain}

\newtheorem{corollary}{Corollary}

\newtheorem{definition}{Definition}

\newtheorem{proposition}{Proposition}

\newtheorem{theorem}{Theorem}
\numberwithin{equation}{section}
\begin{document}
\title[$\phi-$ {Classical 1-Absorbing Prime Submodules}]{$\phi-$ {Classical 1-Absorbing Prime Submodules}}
\author{Zeynep Y\i lmaz U\c{c}ar}
\address{Department of Mathematics, Yildiz Technical University, Istanbul, Turkey;
Orcid: 0009-0009-8107-504X}
\email{zeynepyilmaz@hotmail.com.tr}
\author{Bayram Ali Ersoy}
\address{Department of Mathematics, Yildiz Technical University, Istanbul, Turkey;
Orcid: 0000-0002-8307-9644}
\email{ersoya@yildiz.edu.tr}
\author{\"{U}nsal Tekir}
\address{Department of Mathematics, Marmara University, Istanbul, Turkey; \\
Orcid: 0000-0003-0739-1449}
\email{utekir@marmara.edu.tr}
\author{Ece Yetkin \c{C}elikel}
\address{Department of Basic Sciences, Hasan Kalyoncu University, Gaziantep, Turkey;\\
Orcid: 0000-0001-6194-656X}
\email{ece.celikel@hku.edu.tr}
\author{Serkan Onar}
\address{Department of Mathematical Engineering, Yildiz Technical University, Istanbul,
Turkey; Orcid: 0000-0003-3084-7694}
\email{serkan10ar@gmail.com}

\begin{abstract}
In this paper, all rings are commutative with nonzero identity. Let $R$ be a
ring and $M$ be an $R$-module. The purpose of this paper is to introduce and
study the class of $\phi-$classical 1-absorbing prime submodules. Let
$\phi:S(M)\longrightarrow S(M)\cup\left\{  \emptyset\right\}  $ be a function
where $S\left(  M\right)  $ is the set of all submodules of $M$. A proper
submodule $N$ of $M$ is called a $\phi-$classical 1-absorbing prime submodule
if whenever nonunits $a,b,c\in R$ and $m\in M$ with $abcm\in N\setminus
\phi(N)$, then $abm\in N$ or $cm\in N$. Several characterizations of $\phi-$
classical 1-absorbing prime submodules are given. Among more, it is verified
that if $\ N_{1}N_{2}N_{3}N_{4}\subseteq N$ for some proper submodules
$N_{1},N_{2},N_{3}$ and for some submodule $N_{4}$ of multiplication module
$M$ \ such that $N_{1}N_{2}N_{3}N_{4}\nsubseteq\phi\left(  N\right)  $, then
either $N_{1}N_{2}N_{4}\subseteq N$ or $N_{3}N_{4}\subseteq N$. 

\end{abstract}
\maketitle

\section{Introduction}

In this paper, all rings are commutative with nonzero identity and all modules
are considered to be unitary. The concept of prime ideals and its
generalizations have a significant place in commutative algebra since they are
used in understanding the structure of rings. In \cite{X9}, in order to study
unique factorization domain, Bhatwadekar and Sharma defined a new class of
prime ideals. A proper ideal $I$ is called almost prime ideal if $ab\in I$
$\backslash$ $I^{2}$ for some $a,b\in R$ implies that either $a\in I$ or $b\in
I$. They investigated the relations among the prime ideals, pseudo prime
ideals and almost prime ideals of $R$. Anderson and Batanieh in \cite{X2} gave
a generalization of prime ideals which covers mentioned definition above. Let
$\phi:J(R)\longrightarrow J(R)\cup\left\{  0\right\}  $ be a function where
$J(R)$ is the set of all ideals of $R$. A proper ideal $I$ of $R$ is said to
be $\phi-$prime if for $a,b\in R$ with $ab\in I$ $\backslash$ $\phi(I)$, then
$a\in I$ or $b\in I$. Recently, in \cite{X13}, Koc et al. defined weakly
1-absorbing prime ideal. A proper ideal $I$ of $R$ is said to be a weakly
1-absorbing prime ideal if whenever nonunits $a,b,c\in R$ with $0\neq abc\in
I$, implies that $ab\in I$ or $c\in I$. Several authors have extended the
notion of prime ideals to modules. See, for example, \cite{X10,X11,X15,X16}.
Let $M$ be a module over a commutative ring $R$. A proper submodule $N$ of $M$
is called prime if for $a\in R$ and $m\in M$, $am\in N$ implies that $m\in N$
or $a\in$ $(N:_{R}M)=\{r\in R\mid rM\subseteq N\}$. Zamani \cite{X23}
introduced the concept of $\phi-$prime submodules. Let $\phi:$
$S(M)\longrightarrow$ $S(M)\cup\left\{  \emptyset\right\}  $ be a function
where $S(M)$ is the set of all submodules of $M$. A proper submodule $N$ of an
$R-$module $M$ is called $\phi-$prime if $a\in R$ and $m\in M$ with $am\in N$
$\ \backslash$ $\phi(N)$, then $m\in N$ or $a\in(N:_{R}M)$. He defined the map
$\phi_{\alpha}:$ $S(M)\longrightarrow S(M)\cup\left\{  \emptyset\right\}  $ as follows:

\begin{enumerate}
\item \ $\phi_{\emptyset}$ $:$ \ $\phi(N)=\emptyset$ defines prime submodules.

\item \ $\phi_{0}$ $:$ \ $\phi(N)=\left\{  0\right\}  $ defines weakly prime submodules.

\item \ $\phi_{2}$ $:$ \ $\phi(N)=(N:_{R}M)N$ defines almost prime submodules.

\item \ $\phi_{n}(n\geq2)$ $:$ \ $\phi(I)=(N:_{R}M)^{n-1}N$ defines $n-$almost
prime submodules.

\item \ $\phi_{\omega}$ $:$ \ $\phi(N)=\cap_{n=1}^{\infty}(N:_{R}M)^{n}N$
defines $\omega-$prime submodules.

\item \ $\phi_{1}$ $:$ \ $\phi(N)=N$ defines any submodules.
\end{enumerate}

Also, Moradi and Azizi \cite{X17} investigated the notion of $n$-almost prime
submodules. A proper submodule $N$ of $M$ is called a classical prime
submodule, if for each $a,b\in R$ and $m\in M$, $abm\in N$ implies that $am\in
N$ or $bm\in N$. This notion of classical prime submodules has been
extensively studied by Behboodi in \cite{X5-X7}. For more information on
classical prime submodules, the reader is referred to \cite{X3,X4,X8,X18,X19}.
As a recent research, classical 1-absorbing prime submodules are presented and
investigated in \cite{U}. A proper submodule $N$ of $M$ is said to be a
classical 1-absorbing prime submodule, if for each $m\in M$ and nonunits
$a,b,c\in R$ with $abcm\in N$ implies that $abm\in N$ or $cm\in N$. 

In this paper, we define and study on $\phi-$ classical 1-absorbing prime
submodules which is a generalization of classical 1-absorbing prime
submodules. Throughout this paper $\phi:$ $S(M)\longrightarrow$ $S(M)\cup
\left\{  \emptyset\right\}  $ denotes a function. Since $N\setminus
\phi(N)=N\setminus\left(  N\cap\phi(N)\right)  $, for any submodule $N$ of
$M$, without loss of generality we may assume that $\phi(N)\subseteq N$. For
any two functions $\psi_{1},\psi_{2}:S(M)\longrightarrow$ $S(M)\cup\left\{
\emptyset\right\}  $, we say $\psi_{1}\leq\psi_{2}$ if $\psi_{1}%
(N)\subseteq\psi_{2}(N)$ for each $N\in S(M)$. Thus clearly we have the
following order : $\phi_{\emptyset}\leq\phi_{0}\leq\phi_{\omega}\leq
....\leq\phi_{n+1}\leq\phi_{n}\leq.....\leq\phi_{2}\leq\phi_{1}.$ Whenever
$\psi_{1}\leq\psi_{2}$, any $\psi_{1}-$classical 1-absorbing prime submodule
is $\psi_{2}-$classical 1-absorbing prime. Among other results in this paper,
we give several characterizations of $\phi-$classical 1-absorbing prime
submodules in general modules over commutative rings (see, Theorems
\ref{theo5},\ref{theo1},\ref{theo7}). Also, we investigate the stability of
$\phi-$classical 1-absorbing prime submodules under homomorphism, in
localization of modules, in cartesian product of modules, in multiplication
modules (See, Theorems \ref{theo9}, \ref{theo13}, \ref{theo18}- \ref{theo20},
Proposition \ref{theo3}). Furthermore, we determine the $\phi-$classical
1-absorbing prime submodules of tensor product $F\otimes M$ for a flat
(faithfully flat) $R-$module $F$ and any $R-$module $M$ (see, Theorem
\ref{theo10}, Corollary \ref{cor4}).

\section{Properties of $\phi-$Classical 1-Absorbing Prime Submodules}

\begin{definition}
Let $M$ be an $R$-module and $\phi:$ $S(M)\longrightarrow$ $S(M)\cup\left\{
\emptyset\right\}  $ be a function where $S(M)$ is the set of all submodules
of $M$. Let $N$ be a proper submodule of $M$. Then we say that $N$ is a
$\phi-$classical 1-absorbing prime submodule of $M$ if whenever nonunits
$a,b,c\in R$ and $m\in M$ with $abcm\in N\setminus\phi(N)$, then $abm\in N$ or
$cm\in N$. We define the map $\phi_{\alpha}:$ $S(M)\longrightarrow
S(M)\cup\left\{  \emptyset\right\}  $ as follows:
\end{definition}

\begin{enumerate}
\item \ $\phi_{\emptyset}$ $:$ \ $\phi(N)=\emptyset$ defines classical
1-absorbing prime submodules.

\item \ $\phi_{0}$ $:$ \ $\phi(N)=\left\{  0\right\}  $ defines weakly
classical 1-absorbing prime submodules.

\item \ $\phi_{2}$ $:$ \ $\phi(N)=(N:_{R}M)N$ defines almost classical
1-absorbing prime submodules.

\item \ $\phi_{n}(n\geq2)$ $:$ \ $\phi(I)=(N:_{R}M)^{n-1}N$ defines $n-$almost
classical 1-absorbing prime submodules.

\item \ $\phi_{\omega}$ $:$ \ $\phi(N)=\cap_{n=1}^{\infty}(N:_{R}M)^{n}N$
defines $\omega-$classical 1-absorbing prime submodules.

\item \ $\phi_{1}$ $:$ \ $\phi(N)=N$ defines any submodules.
\end{enumerate}

Let $R$ be a ring and $\psi:J(R)\longrightarrow J(R)\cup\left\{
\emptyset\right\}  $ be a function. Clearly, $I$ is a $\psi-$1-absorbing prime
ideal of $R$ if and only if $I$ is a $\psi-$classical 1-absorbing prime
submodule of $R$-module $R$. 

In \cite{X20}, Quartararo et al. said that a commutative ring $R$ is a
$u$-ring provided $R$ has the property that an ideal contained in a finite
union of ideals must be contained in one of those ideals; and $um-$ring is a
ring $R$ with the property that an $R$-module which is equal to a finite union
of submodules must be equal one of them. They show that every B\'{e}zout ring
is a $u-$ring. Moreover, they proved that every Pr\"{u}fer domain is a
$u-$domain. Also, any ring which contains an infinite field as a subring is a
$u-$ring \cite[Exercise 3.63]{X21}.

\begin{theorem}
\label{theo4}Let $M$ be an $R-$module and $\phi:S(M)\longrightarrow
S(M)\cup\left\{  \emptyset\right\}  $ be a function. Suppose that $N$ is a
$\phi-$classical 1-absorbing prime submodule of $M$. Then
\end{theorem}

\begin{enumerate}
\item For every nonunits $a,b,c\in R$ and $m\in M$, $\left(  N:_{R}%
abcm\right)  =\left(  \phi\left(  N\right)  :_{R}abcm\right)  \cup\left(
N:_{R}abm\right)  \cup\left(  N:_{R}cm\right)  $.

\item If $R$ is a $u-$ring, then for every nonunits $a,b,c\in R$ and $m\in M$,
$\left(  N:_{R}abcm\right)  =\left(  \phi\left(  N\right)  :_{R}abcm\right)  $
or $\left(  N:_{R}abcm\right)  =\left(  N:_{R}abm\right)  $ or $\left(
N:_{R}abcm\right)  =\left(  N:_{R}cm\right)  $.
\end{enumerate}

\begin{proof}
$\left(  1\right)  $ Let nonunits $a,b,c\in R$ and $m\in M$. Suppose that
$r\in\left(  N:_{R}abcm\right)  $. Then $abc\left(  rm\right)  \in N$. If
$abc\left(  rm\right)  \in\phi\left(  N\right)  $, then $r\in\left(
\phi\left(  N\right)  :_{R}abcm\right)  $. Now, assume that $abc\left(
rm\right)  \notin\phi\left(  N\right)  $. Then we have $abc\left(  rm\right)
\in N\setminus\phi(N)$ and since $N$ is a $\phi-$classical 1-absorbing prime
submodule we get $ab\left(  rm\right)  \in N$ or $c\left(  rm\right)  \in N$.
Hence $r\in\left(  N:_{R}abm\right)  $ or $r\in\left(  N:_{R}cm\right)  $.
Thus, we conclude that $r\in\left(  \phi\left(  N\right)  :_{R}abcm\right)
\cup\left(  N:_{R}abm\right)  \cup\left(  N:_{R}cm\right)  $, that is $\left(
N:_{R}abcm\right)  \subseteq\left(  \phi\left(  N\right)  :_{R}abcm\right)
\cup\left(  N:_{R}abm\right)  \cup\left(  N:_{R}cm\right)  $. Since the
reverse inclusion is always true, we have the required equality.

$\left(  2\right)  $ Since $R$ is a $u$-ring, apply part (1).
\end{proof}

Let $M$ be an $R-$module and $N$ a submodule of $M$. For every $a\in R$,
$\left\{  m\in M\mid am\in N\right\}  $ is denoted by $\left(  N:_{M}a\right)
$. It is easy to see that $\left(  N:_{M}a\right)  $ is a submodule of $M$
containing $N$. In the next theorem, we obtain several equivalent statements
to characterize $\phi-$classical 1-absorbing prime submodules.

\begin{theorem}
\label{theo5}Let $M$ be an $R$-module, $\phi:S(M)\longrightarrow
S(M)\cup\left\{  \emptyset\right\}  $ a function and $N$ be a proper submodule
$M$. The following conditions are equivalent.
\end{theorem}

\begin{enumerate}
\item $N$ is a $\phi-$classical 1-absorbing prime submodule of $M$.

\item For every nonunits $a,b,c\in R$, $\left(  N:_{M}abc\right)  =\left(
\phi\left(  N\right)  :_{M}abc\right)  \cup\left(  N:_{M}ab\right)
\cup\left(  N:_{M}c\right)  $.

\item For every nonunits $a,b\in R$ and $m\in M$ with $abm\notin N$; $\left(
N:_{R}abm\right)  =\left(  \phi\left(  N\right)  :_{R}abm\right)  \cup\left(
N:_{R}m\right)  $.

\item For every nonunits $a,b\in R$ and $m\in M$ with $abm\notin N$; $\left(
N:_{R}abm\right)  =\left(  \phi\left(  N\right)  :_{R}abm\right)  $ or
$\left(  N:_{R}abm\right)  =\left(  N:_{R}m\right)  $.

\item For every nonunits $a,b\in R$, every ideal $I$ of $R$ and $m\in M$ with
$abIm\subseteq N$ and $abIm\nsubseteq\phi\left(  N\right)  $, either $abm\in
N$ or $Im\subseteq N$.

\item For every nonunits $a,b\in R$, $m\in M$ and every proper ideal $I$ of
$R$ with $abIm\subseteq N$ and $abIm\nsubseteq\phi\left(  N\right)  $, either
$aIm\subseteq N$ or $bm\in N$.

\item For every proper ideal $I$ of $R$, nonunit $a\in R$ and $m\in M$ with
$aIm\nsubseteq N$, $\left(  N:_{R}aIm\right)  =\left(  \phi\left(  N\right)
:_{R}aIm\right)  $ or $\left(  N:_{R}aIm\right)  =\left(  N:_{R}m\right)  $.

\item For every proper ideals $I,J,K$ of $R$ and $m\in M$ with $IJKm\subseteq
N$ and $IJKm\nsubseteq\phi\left(  N\right)  $, either $IJm\subseteq N$ or
$Km\subseteq N$.
\end{enumerate}

\begin{proof}
$(1)\Rightarrow(2)$ Suppose that $N$ is a $\phi-$classical 1-absorbing prime
submodule of $M$. Let $m\in$ $\left(  N:_{M}abc\right)  $. Then $abcm\in N$.
If $abcm\in\phi\left(  N\right)  $, then $m\in\left(  \phi\left(  N\right)
:_{M}abc\right)  $. Assume that $abcm\notin\phi\left(  N\right)  $. Since $N$
is $\phi-$ classical 1-absorbing prime, we have either $abm\in N$ or $cm\in
N$. Hence $m\in\left(  N:_{M}ab\right)  $ or $m\in\left(  N:_{M}c\right)  $.
Thus we conclude $\left(  N:_{M}abc\right)  \subseteq\left(  \phi\left(
N\right)  :_{M}abc\right)  \cup\left(  N:_{M}ab\right)  \cup\left(
N:_{M}c\right)  $. Since the reverse inclusion is always true, we have the
required equality.

$(2)\Rightarrow(3)$ Suppose that $abm\notin N$ for some nonunits $a,b\in R$
and $m\in M$. Let $x\in\left(  N:_{R}abm\right)  $. Then $m\in\left(
N:_{M}abx\right)  $ and since $abm\notin N$, then $m\notin\left(
N:_{M}ab\right)  $. Thus by part (2) we have $m\in\left(  \phi\left(
N\right)  :_{M}abx\right)  $ or $m\in\left(  N:_{M}x\right)  $ whence
$x\in\left(  \phi\left(  N\right)  :_{R}abm\right)  $ or $x\in\left(
N:_{R}m\right)  $. Then we get $\left(  N:_{R}abm\right)  =\left(  \phi\left(
N\right)  :_{R}abm\right)  \cup\left(  N:_{R}m\right)  $.

$(3)\Longrightarrow(4)$ By the fact that if an ideal (a subgroup) is the union
of two ideals (two subgroups), then it is equal to one of them.

$(4)\Longrightarrow(5)$ Let $abIm\subseteq N$ and $abIm\nsubseteq\phi\left(
N\right)  $ for some nonunits $a,b\in R$, an ideal $I$ of $R$ and $m\in M$.
Hence $I\subseteq(N:_{R}abm)$ and $I\nsubseteq(\phi\left(  N\right)
:_{R}abm)$. If $abm\in N$, then we are done. So, assume that $abm\notin N$.
Therefore by part (4) we have that $I\subseteq\left(  N:_{R}m\right)  $, and
then $Im\subseteq N$.

$(5)\Longrightarrow(6)$ Let $abIm\subseteq N$ and $abIm\nsubseteq\phi\left(
N\right)  $ for some proper ideal $I$ of $R$, nonunits $a,b\in R$ and $m\in
M$. Assume that $aIm\nsubseteq N$. Then there exists $x\in I$ such that
$axm\notin N$. Then note that $ax\left(  Rb\right)  m\subseteq N$. If
$ax\left(  Rb\right)  m\nsubseteq\phi\left(  N\right)  $, then by part (5) we
have $(Rb)m\subseteq N$ which completes the proof. So assume that $ax\left(
Rb\right)  m\subseteq\phi\left(  N\right)  $. Since $abIm\nsubseteq\phi\left(
N\right)  $, there exists $y\in I$ such that $ay\left(  Rb\right)
m\nsubseteq\phi\left(  N\right)  $. This implies that $a\left(  x+y\right)
\left(  Rb\right)  m\subseteq N\setminus\phi\left(  N\right)  $. As $ay\left(
Rb\right)  m\subseteq N\setminus\phi\left(  N\right)  $, again by part (5), we
have $aym\in N$ or $\left(  Rb\right)  m\subseteq N$. Let $aym\in N$. Then we
have $a\left(  x+y\right)  m\notin N$. Since $I$ is a proper ideal, we
conclude that $x+y$ is a nonunit. As $a\left(  x+y\right)  \left(  Rb\right)
m\subseteq N\setminus\phi\left(  N\right)  $, by part (5), we get
$(Rb)m\subseteq N$ which implies that $bm\in N$. In other case, we have $bm\in
N$ which completes the proof.

$(6)\Longrightarrow(7)$ Let $aIm\nsubseteq N$ for some nonunit $a\in R$, some
proper ideal $I$ of $R$ and $m\in M$. Choose $b\in$ $\left(  N:_{R}aIm\right)
$. Then we have $abIm\subseteq N$. If $abIm\subseteq\phi\left(  N\right)  $,
then we get $b\in\left(  \phi\left(  N\right)  :_{R}aIm\right)  $. Now assume
that $abIm\nsubseteq\phi\left(  N\right)  $. Then by part (6), we get $bm\in
N$ which implies that $b\in\left(  N:_{R}m\right)  $. Thus, we conclude that
$\left(  N:_{R}aIm\right)  =\left(  \phi\left(  N\right)  :_{R}aIm\right)
\cup\left(  N:_{R}m\right)  $. This implies that $\left(  N:_{R}aIm\right)
=\left(  \phi\left(  N\right)  :_{R}aIm\right)  $ or $\left(  N:_{R}%
aIm\right)  =\left(  N:_{R}m\right)  $.

$(7)\Longrightarrow(8)$ Let $IJKm\subseteq N$ and $IJKm\nsubseteq\phi\left(
N\right)  $ for some proper ideals $I,J,K$ of $R$ and $m\in M$. Assume that
$IJm\nsubseteq N$. Then there exists $a\in J$ such that $aIm\nsubseteq N$.
Since $aIKm\subseteq N$, we have $K\subseteq\left(  N:_{R}aIm\right)  =\left(
\phi\left(  N\right)  :_{R}aIm\right)  $ or $K\subseteq\left(  N:_{R}%
aIm\right)  =\left(  N:_{R}m\right)  $ by part (7). This gives $aIKm\subseteq
\phi\left(  N\right)  $ or $Km\subseteq N$. In the latter case, we are done.
So, we may assume that $aIKm\subseteq\phi\left(  N\right)  $. Since
$IJKm\nsubseteq\phi\left(  N\right)  $, there exists $b\in J$ such that
$IbKm\nsubseteq\phi\left(  N\right)  $. Since $K\subseteq\left(
N:_{R}bIm\right)  $, again by part (7), we have $bIm\subseteq N$ or
$Km\subseteq N$. Now assume that $bIm\subseteq N$. Then we have $\left(
a+b\right)  Im\nsubseteq N$. Since $J$ is a proper ideal, our conclusion is
that $a+b$ is a nonunit. As $\left(  a+b\right)  IKm\subseteq N$, $\left(
a+b\right)  IKm\nsubseteq\phi\left(  N\right)  $ we obtain that $K\subseteq
\left(  N:_{R}\left(  a+b\right)  Im\right)  $ and by part (7), we conclude
that $Km\subseteq N$.

$(8)\Longrightarrow(1)$ Clear.
\end{proof}

Let $M$ be an $R-$module, $K$ a submodule of $M$ and $\phi:$
$S(M)\longrightarrow$ $S(M)\cup\left\{  \emptyset\right\}  $ a function.
Define $\phi_{K}:$ $S(M/K)\longrightarrow$ $S(M/K)\cup\left\{  0\right\}  $ by
$\phi_{K}(N/K)=\left(  \phi(N)+K\right)  /K$ for every $N\in S(M)$ with
$N\supseteq K$ (and $\phi_{K}(N/K)=\emptyset$ if $\phi(N)=\emptyset$).

\begin{theorem}
$\label{theo1}$Let $M$ be an $R$-module and $K\subseteq N$ be proper
submodules of $M$. Suppose that $\phi:$ $S(M)\longrightarrow$ $S(M)\cup
\left\{  \emptyset\right\}  $ be a function.
\end{theorem}

\begin{enumerate}
\item If $N$ is a $\phi-$classical 1-absorbing prime submodule of $M$, then
$N/K$ is a $\phi_{K}-$classical 1-absorbing prime submodule of $M/K$.

\item If $K\subseteq\phi(N)$ and $N/K$ is a $\phi_{K}-$classical 1-absorbing
prime submodule of $M/K$, then $N$ is a $\phi-$classical 1-absorbing prime
submodule of $M$.

\item If $\phi(N)\subseteq K$ and $N$ is a $\phi-$classical 1-absorbing prime
submodule of$M$, then $N/K$ is a weakly classical 1-absorbing prime submodule
of $M/K$.

\item If $\phi(K)\subseteq\phi(N)$, $K$ is a $\phi-$classical 1-absorbing
prime submodule of $M$ and $N/K$ is a weakly classical 1-absorbing prime
submodule of $M/K$, then $N$ is a $\phi-$classical 1-absorbing prime submodule
of $M$.

\item If $N$ is a $\phi-$classical 1-absorbing prime submodule of $M$ and
$\phi\left(  N\right)  $ be a classical 1-absorbing prime submodule of $M$,
then $N$ is a classical 1-absorbing prime submodule of $M$.
\end{enumerate}

\begin{proof}
$(1)$ Let nonunits $a,b,c\in R$ and $m\in M$ such that $abc(m+K)\in
N/K\setminus\phi_{K}(N/K)$. Then $abcm\in N\setminus\phi(N)$. Since $N$ is a
$\phi-$classical 1-absorbing prime submodule of $M$ we have either $abm\in N$
or $cm\in N$. Therefore, $ab(m+K)\in N/K$ or $c(m+K)\in N/K$. Consequently,
$N/K$ is a $\phi_{K}-$classical 1-absorbing prime submodule of $M/K$.

$\left(  2\right)  $ Let nonunits $a,b,c\in R$, $\phi_{K}\left(  N/K\right)
=\left(  \phi\left(  N\right)  /K\right)  $ and $m\in M$ such that $abcm\in
N\setminus\phi(N)$. Then $abc(m+K)\in N/K\setminus\left(  \phi\left(
N\right)  /K\right)  $. Since $N/K$ is a $\phi_{K}-$classical 1-absorbing
prime submodule we get either $ab\left(  m+K\right)  \in N/K$ or $c\left(
m+K\right)  \in N/K$, and so $abm\in N$ or $cm\in N$.

$\left(  3\right)  $ Let nonunits $a,b,c\in R$ and $m\in M$ such that $0\neq
abc\left(  m+K\right)  \in N/K$. Hence $abcm\in N\setminus\phi\left(
N\right)  $ as $\phi\left(  N\right)  \subseteq K$. Since $N$ is a $\phi
-$classical 1-absorbing prime submodule we have either $abm\in N$ or $cm\in
N$. Therefore, $ab\left(  m+K\right)  \in N/K$ or $c\left(  m+K\right)  \in
N/K$. Thus, $N/K$ is a weakly classical 1-absorbing prime submodule of $M/K$.

$\left(  4\right)  $ Let nonunits $a,b,c\in R$ and $m\in M$ such that $abcm\in
N\setminus\phi(N)$. Note that $\phi\left(  K\right)  \subseteq\phi\left(
N\right)  $ implies that $abcm\notin\phi\left(  K\right)  $. If $abcm\in K$,
then we have $abcm\in K\setminus\phi\left(  K\right)  $. Since $K$ is a
$\phi-$classical 1-absorbing prime submodule we have either $abm\in K\subseteq
N$ or $cm\in K\subseteq N$. Now, assume that $abcm\notin K$. So $0\neq
abc\left(  m+K\right)  \in N/K$. Since $N/K$ is a weakly classical 1-absorbing
prime submodule of $M/K$ we have either $ab\left(  m+K\right)  \in N/K$ or
$c\left(  m+K\right)  \in N/K$. Thus $abm\in N$ or $cm\in N$. Consequently,
$N$ is a $\phi-$classical 1-absorbing prime submodule of $M$.

$\left(  5\right)  $ Let $N$ be a $\phi-$classical 1-absorbing prime submodule
of $M$. Assume that $abcm\in N$ for some nonunits $a,b,c\in R$ and $m\in M$.
If $abcm\in\phi\left(  N\right)  $, then since $\phi\left(  N\right)  $ is
classical 1-absorbing prime, we conclude that $abm\in\phi\left(  N\right)
\subseteq N$ or $cm\in\phi\left(  N\right)  \subseteq N$, and so we are done.
When $abcm\notin\phi\left(  N\right)  $ clearly the result follows.
\end{proof}

In view of \ref{theo1} (2) and (3), we conclude that a proper submodule $N$ of
$M$ is a $\phi-$classical 1-absorbing prime submodule if and only if
$N/\phi\left(  N\right)  $ is a weakly classical 1-absorbing prime submodule
of $M/$ $\phi\left(  N\right)  $.

\begin{theorem}
$\label{theo2}$ Let $M$ be an $R-$module and $N$ be a proper submodule of $M$.
Suppose that $\phi:S(M)\longrightarrow S(M)\cup\left\{  \emptyset\right\}  $
and $\psi:J(R)\longrightarrow J(R)\cup\left\{  \emptyset\right\}  $ be two functions.
\end{theorem}

\begin{enumerate}
\item If $N$ is a $\phi-$classical 1-absorbing prime submodule of $M$, then
$(N:_{R}m)$ is a $\psi-$1-absorbing prime ideal of $R$ for every $m\in
M\setminus N$ satisfying $\left(  \phi\left(  N\right)  :_{R}m\right)
\subseteq\psi\left(  (N:_{R}m)\right)  $.

\item If $(N:_{R}m)$ is a $\psi-$1-absorbing prime ideal of $R$ for every
$m\in M\setminus N$ with $\psi\left(  (N:_{R}m)\right)  \subseteq\left(
\phi\left(  N\right)  :_{R}m\right)  $, then $N$ is a $\phi-$classical
1-absorbing prime submodule of $M$.
\end{enumerate}

\begin{proof}
$\left(  1\right)  $ Suppose that $N$ is a $\phi-$classical 1-absorbing prime
submodule. Let $m\in M\setminus N$ with $\left(  \phi\left(  N\right)
:_{R}m\right)  \subseteq\psi\left(  (N:_{R}m)\right)  $ and $abc\in
(N:_{R}m)\setminus\psi\left(  (N:_{R}m)\right)  $ for some nonunits $a,b,c\in
R$. Then $abcm\in N\setminus\phi\left(  N\right)  $. Since $N$ is a $\phi
-$classical 1-absorbing prime submodule we have either $abm\in N$ or $cm\in
N$. Hence $ab\in(N:_{R}m)$ or $c\in(N:_{R}m)$. Consequently, $(N:_{R}m)$ is a
$\psi-$1-absorbing prime ideal of $R$.

$\left(  2\right)  $ Assume that $(N:_{R}m)$ is a $\psi-$1-absorbing prime
ideal of $R$ for every $m\in M\setminus N$ with $\psi\left(  (N:_{R}m)\right)
\subseteq\left(  \phi\left(  N\right)  :_{R}m\right)  $. Let $abcm\in
N\setminus\phi(N)$ for some $m\in M$ and nonunits $a,b,c\in R$. If $m\in N$,
then we are done. So we assume that $m\notin N$. Hence $abc\in(N:_{R}%
m)\setminus\left(  \phi\left(  N\right)  :_{R}m\right)  $. Since $\psi\left(
(N:_{R}m)\right)  \subseteq\left(  \phi\left(  N\right)  :_{R}m\right)  $ we
have $abc\in(N:_{R}m)\setminus\psi\left(  (N:_{R}m)\right)  $. Since
$(N:_{R}m)$ is a $\psi-$1-absorbing prime ideal we have either $ab\in
(N:_{R}m)$ or $c\in(N:_{R}m)$. Hence $abm\in N$ or $cm\in N$ and so $N$ is a
$\phi-$classical 1-absorbing prime submodule of $M$.
\end{proof}

E.A.U\u{g}urlu \cite{X22}, generalized the concept of prime submodules of a
module over a commutative ring as follows: Let $N$ be a proper submodule of an
$R-$module $M$. Then $N$ is said to be a 1-absorbing prime submodule of $M$ if
whenever nonunits $a,b\in R$ and $m\in M$ with $abm\in N$, then $m\in N$ or
$ab\in(N:_{R}M).$ We recall that an $R$-module $M$ is called a multiplication
module if every submodule $N$ of $M$ has the form $IM$ for some ideal $I$ of
$R$, \cite{X12}. Note that, since $I\subseteq(N:_{R}M)$, $N=IM\subseteq
(N:_{R}M)M\subseteq N$, and $N=(N:_{R}M)M$. Let $N$ and $K$ be submodules of a
multiplication $R-$module $M$ with $N=I_{1}M$ and $K=I_{2}M$ for some ideals
$I_{1}$ and $I_{2}$ of $R$. The product of $N$ and $K$ denoted by $NK$ is
defined by $NK=I_{1}I_{2}M$. Then by \cite[Theorem 3.4]{X1}, the product of
$N$ and $K$ is independent of presentations of $N$ and $K$. Let $N$ be a
proper submodule of a nonzero $R$-module $M$. Then $M$-radical of $N$, denoted
by $M$-$rad(N)$, is defined to be intersection of all prime submodules of $M$
containing $N$. If $M$ has no prime submodule containing $N$,then
$M$-rad$(N)=M$. It is shown in \cite[Theorem 2.12]{X12} that if $N$ is a
proper submodule of a multiplication $R$-module $M$, then $M$-rad$(N)=\sqrt
{(N:_{R}M)}M$. 

\begin{definition}
Let $M$ be an $R-$module, $\phi:S(M)\longrightarrow S(M)\cup\left\{
\emptyset\right\}  $ a function where $S\left(  M\right)  $ is the set of all
submodules of $M$. Let $N$ be a proper submodule of $M$. $N$ is called $\phi
-$1-absorbing prime submodule of $M$ if whenever nonunits $a,b\in R$ and $m\in
M$ with $abm\in N$ and $abm\notin\phi\left(  N\right)  $, then $m\in N$ or
$ab\in(N:_{R}M)$. We can define the following special functions $\phi_{\alpha
}$ as follows. Let $N$ be a $\phi_{\alpha}-$1-absorbing prime submodule of a
multiplication $R-$module $M$. Then
\end{definition}

$\phi_{\emptyset}$ $(N)=\emptyset$ 1-absorbing prime submodules.

$\phi_{0}$ $(N)=0$ weakly 1-absorbing prime submodules.

$\phi_{2}$ $(N)=N^{2}$ almost 1-absorbing prime submodules.

\ \ \ $...$

$\phi_{n}(N)=N^{n}$ $n$-almost 1-absorbing prime submodules.

$\phi_{\omega}$ $(N)=\cap_{n=1}^{\infty}N^{n}$ $\omega-$1-absorbing prime
submodules.\medskip

It is clear that if $N$ is a $\phi-$1-absorbing prime submodule of $M$, then
$N$ is a $\phi-$classical 1-absorbing prime submodule of $M.$

\begin{definition}
Let $N$ be a proper submodule of a multiplication $R-$module $M$ and $n\geq2$.
Then $N$ is said to be $n-$potent classical 1-absorbing prime if whenever
nonunits $a,b,c\in R$ and $m\in M$ with $abcm\in N^{n}$, then $abm\in N$ or
$cm\in N$.
\end{definition}

\begin{proposition}
\label{prop2}Let $M$ be a multiplication $R-$module. If $N$ is an $n-$almost
classical 1-absorbing prime submodule of $M$ for some $n\geq2$ and $N$ is a
$k-$potent classical 1-absorbing prime for some $k\leq n$, then $N$ is a
classical 1-absorbing prime submodule of $M$.
\end{proposition}

\begin{proof}
Assume that $N$ is an $n-$almost classical 1-absorbing prime submodule of $M$.
Let $abcm\in N$ for some nonunits $a,b,c\in R$ and $m\in M$. If $abcm\notin
N^{k}$, then $abcm\notin N^{n}$. Since $N$ is an $n-$almost classical
1-absorbing prime we have $abm\in N$ or $cm\in N$. Now assume that $abcm\in
N^{k}$. Since $N$ is a $k-$potent classical 1-absorbing prime we get $abm\in
N$ or $cm\in N$. Consequently, $N$ is a classical 1-absorbing prime submodule
of $M$.
\end{proof}

\begin{proposition}
\label{prop3}Let $M$ be a cyclic $R-$module and $\phi:S(M)\longrightarrow
S(M)\cup\left\{  \emptyset\right\}  $ be a function. Then a proper submodule
$N$ of $M$ is a $\phi-$1-absorbing prime submodule if and only if it is a
$\phi-$classical 1-absorbing prime submodule of $M$.
\end{proposition}

\begin{proof}
Let $M=Rm$ for some $m\in M$ and $N$ be a $\phi-$classical 1-absorbing prime
submodule of $M$. Suppose that $abx\in N\setminus\phi\left(  N\right)  $ for
some nonunits $a,b\in R$ and $x\in M$. Then there exists a nonunit element
$c\in R$ such that $x=cm$. Hence $abcm\in N\setminus\phi(N)$ implies either
$abm\in N$ or $cm\in N$. Hence $ab\in(N:_{R}M)$ or $x\in N$. Consequently, $N$
is a $\phi-$1-absorbing prime submodule of $M$. The converse part is straighforward.
\end{proof}

\begin{proposition}
$\label{theo3}$Let $f\colon M\longrightarrow M%
\acute{}%
$ be an epimorphism of $R$-modules and $\phi:S(M)\longrightarrow
S(M)\cup\left\{  \emptyset\right\}  $ and $\phi%
\acute{}%
:S(M%
\acute{}%
)\longrightarrow S(M%
\acute{}%
)\cup\left\{  \emptyset\right\}  $ be two functions. Then the following
conditions hold:
\end{proposition}

\begin{enumerate}
\item If \ $N%
\acute{}%
$ is a $\phi%
\acute{}%
-$classical 1-absorbing prime submodule of $M%
\acute{}%
$ and $\phi\left(  f^{-1}\left(  N%
\acute{}%
\right)  \right)  =f^{-1}\left(  \phi%
\acute{}%
\left(  N%
\acute{}%
\right)  \right)  $, then $f^{-1}\left(  N%
\acute{}%
\right)  $ is a $\phi-$classical 1-absorbing prime submodule of $M$.

\item If $N$ is a $\phi-$classical 1-absorbing prime submodule of \ $M$
containing $Ker(f)$ and $\phi%
\acute{}%
\left(  f\left(  N\right)  \right)  =f\left(  \phi\left(  N\right)  \right)
$, then $f\left(  N\right)  $ is a $\phi%
\acute{}%
-$classical 1-absorbing prime submodule of $M%
\acute{}%
$.
\end{enumerate}

\begin{proof}
$\left(  \text{1}\right)  $ Since $f$ is an epimorphism, $f^{-1}\left(  N%
\acute{}%
\right)  $ is a proper submodule of $M$. Let nonunits $a,b,c\in R$ and $m\in
M$ such that $abcm\in f^{-1}\left(  N%
\acute{}%
\right)  \setminus\phi(f^{-1}\left(  N%
\acute{}%
\right)  )$. Since $abcm\in f^{-1}\left(  N%
\acute{}%
\right)  $, we have $f(abcm)=abcf(m)\in N%
\acute{}%
$. Also, $\phi\left(  f^{-1}\left(  N%
\acute{}%
\right)  \right)  =f^{-1}\left(  \phi%
\acute{}%
\left(  N%
\acute{}%
\right)  \right)  $ implies that $abcm\notin\phi%
\acute{}%
\left(  N%
\acute{}%
\right)  $. Thus, $abcf\left(  m\right)  \in N%
\acute{}%
\setminus\phi%
\acute{}%
\left(  N%
\acute{}%
\right)  $ implies either $abf\left(  m\right)  \in N%
\acute{}%
$ or $cf\left(  m\right)  \in N%
\acute{}%
$. Hence $abm\in f^{-1}\left(  N%
\acute{}%
\right)  $ or $cm\in f^{-1}\left(  N%
\acute{}%
\right)  $ as needed.

$\left(  2\right)  $ Let nonunits $a,b,c\in R$ and $m%
\acute{}%
\in M%
\acute{}%
$ such that $abcm%
\acute{}%
\in f\left(  N\right)  \setminus\phi%
\acute{}%
\left(  f\left(  N\right)  \right)  .$ Since $f$ $\ $is epimorphism, there
exists $m\in M$ such that $m%
\acute{}%
=f\left(  m\right)  $. Then $abcf\left(  m\right)  \in f\left(  N\right)
\setminus\phi%
\acute{}%
\left(  f\left(  N\right)  \right)  $. Hence $f\left(  abcm\right)  \in
f\left(  N\right)  \setminus\phi%
\acute{}%
\left(  f\left(  N\right)  \right)  $. Since $Ker\left(  f\right)  \subseteq
N$, we have $abcm\in N$. Since $\phi%
\acute{}%
\left(  f\left(  N\right)  \right)  =f\left(  \phi\left(  N\right)  \right)  $
we get $abcm\notin\phi\left(  N\right)  $. Thus $\ abcm\in N\setminus\phi(N)$
and $N$ is a $\phi-$ classical 1-absorbing prime submodule of \ $M$ we have
either $abm\in N$ or $cm\in N$. Then $abf\left(  m\right)  \in f\left(
N\right)  $ or $cf\left(  m\right)  \in f\left(  N\right)  $ and so $abm%
\acute{}%
\in f\left(  N\right)  $ or $cm%
\acute{}%
\in f\left(  N\right)  $. Consequently, $f\left(  N\right)  $ is a $\phi%
\acute{}%
-$classical 1-absorbing prime submodule of $M%
\acute{}%
$.
\end{proof}

Let $M$ be a multiplication $R-$module and $K,L$ be submodules of $M$. Then
there are ideals $I,J$ of $R$ such that $K=IM$ and $L=JM$. Thus $KL=IJM=IL$.
In particular $KM=IM=K$. Also, for any $m\in M$ we define $Km:=KRm$. Hence
$Km=IRm=Im$. Next, we characterize $\phi-$ classical 1-absorbing prime
submodule of multiplication modules in terms of some submodules of them.

\begin{theorem}
\label{theo6} Let $M$ be a multiplication $R$-module, $N$ be a proper
submodule of $M$ and $\phi:S(M)\longrightarrow S(M)\cup\left\{  \emptyset
\right\}  $ be a function. Then the following conditions are equivalent:

\begin{enumerate}
\item $N$ is a $\phi-$ classical 1-absorbing prime submodule of $M$.

\item If$\ N_{1}N_{2}N_{3}m\subseteq N$ for some submodules $N_{1},N_{2}%
,N_{3}$ of $M$ and $m\in M$ such that $N_{1}N_{2}N_{3}m\nsubseteq\phi\left(
N\right)  $, then either $N_{1}N_{2}m\subseteq N$ or $N_{3}m\subseteq N$ .
\end{enumerate}
\end{theorem}

\begin{proof}
$(1)\Longrightarrow(2)$ Let $N_{1}N_{2}N_{3}m\subseteq N$ for some submodules
$N_{1},N_{2},N_{3}$ \ of $M$ and $m\in M$ . Since $M$ is a multiplication
$R$-module, there are ideals $I_{1},I_{2},I_{3}$ of $R$ such that $N_{1}%
=I_{1}M$, $N_{2}=I_{2}M$ and $N_{3}=I_{3}M$. Then $N_{1}N_{2}N_{3}m=I_{1}%
I_{2}I_{3}m\subseteq N$ and $I_{1}I_{2}I_{3}m\nsubseteq\phi\left(  N\right)
$. Since by Theorem \ref{theo5} we have either $I_{1}I_{2}m\subseteq N$ or
$I_{3}m\subseteq N$. Therefore $N_{1}N_{2}m\subseteq N$ or $N_{3}m\subseteq N$ .

$(2)\Rightarrow(1)$ Suppose that $I_{1}I_{2}I_{3}m\subseteq N$ for some proper
ideals $I_{1},I_{2},I_{3}$ of $R$ and some $m\in M$. It is sufficient to set
$N_{1}:=I_{1}M$, $N_{2}:=I_{2}M$ and $N_{3}:=I_{3}M$ in part (2).
\end{proof}

Provided that $R$ is a um-ring, we conclude several equivalent conditions.

\begin{theorem}
\label{theo7}Let $R$ be a um-ring, $M$ be an $R$-module and $N$ be a proper
submodule of $M$. Suppose that $\phi:S(M)\longrightarrow S(M)\cup\left\{
\emptyset\right\}  $ be a function. The following conditions are equivalent.

\begin{enumerate}
\item $N$ is a $\phi-$classical 1-absorbing prime submodule of $M$.

\item For every nonunits $a,b,c\in R$, $\left(  N:_{M}abc\right)  =\left(
\phi\left(  N\right)  :_{M}abc\right)  $ or $\left(  N:_{M}abc\right)
=\left(  N:_{M}ab\right)  $ or $\left(  N:_{M}abc\right)  =\left(
N:_{M}c\right)  $.

\item For every nonunits $a,b,c\in R$ and every submodule $L$ of $M$ ;
$abcL\subseteq N$ and $abcL\nsubseteq\phi\left(  N\right)  $ implies that
$abL\subseteq N$ or $cL\subseteq N$.

\item For every nonunits $a,b\in R$ and every submodule $L$ of $M$ with
$abL\nsubseteq N$ , $\left(  N:_{R}abL\right)  =\left(  \phi\left(  N\right)
:_{R}abL\right)  $ or $\left(  N:_{R}abL\right)  =\left(  N:_{R}L\right)  $.

\item For every nonunits $a,b\in R$ , every ideal $J$ of $R$ and every
submodule $L$ of $M$ with $abJL\subseteq N$ and $abJL\nsubseteq\phi\left(
N\right)  $ implies that $abL\subseteq N$ or $JL\subseteq N$.

\item For every nonunit $a\in R,$ every proper ideals $I,J$ of $R$ and every
submodule $L$ of $M$ with $aIJL\subseteq N$ and $aIJL\nsubseteq\phi\left(
N\right)  $ implies that $aIL\subseteq N$ or $JL\subseteq N$.

\item For every proper ideals $I,J$ of $R$ and every submodule $L$ of $M$ with
$IJL\nsubseteq N$ ;$\left(  N:_{R}IJL\right)  =\left(  \phi\left(  N\right)
:_{R}IJL\right)  $ or $\left(  N:_{R}IJL\right)  =\left(  N:_{R}L\right)  $.

\item For every proper ideals $I,J,K$ of $R$ and every submodule $L$ of $M$
with $IJKL\subseteq N$ and $IJKL\nsubseteq\phi\left(  N\right)  $ implies that
$IJL\subseteq N$ or $KL\subseteq N$.
\end{enumerate}
\end{theorem}

\begin{proof}
$(1)\Rightarrow(2)$ We get the claim by Theorem \ref{theo5} and by the
assumption that $R$ is a $um-$ring.

$(2)\Rightarrow(3)$ Let $abcL\subseteq N$ and $abcL\nsubseteq\phi\left(
N\right)  $ for some nonunits $a,b,c\in R$ and submodule $L$ of $M$. Hence
$L\subseteq$ $\left(  N:_{M}abc\right)  $ and $L\nsubseteq\left(  \phi\left(
N\right)  :_{M}abc\right)  $. Therefore, part (2) implies that $\left(
N:_{M}abc\right)  =\left(  N:_{M}ab\right)  $ or $\left(  N:_{M}abc\right)
=\left(  N:_{M}c\right)  $. So, either $L\subseteq\left(  N:_{M}ab\right)  $
or $L\subseteq\left(  N:_{M}c\right)  $, and then $abL\subseteq N$ or
$cL\subseteq N$.

$(3)\Longrightarrow(4)$ Let $abL\nsubseteq N$ for some nonunits $a,b\in R$ and
submodule $L$ of $M$. Assume that $c\in\left(  N:_{R}abL\right)  $. Then $c$
is nonunit and\ $abcL\subseteq N$. If $abcL\subseteq\phi\left(  N\right)  $,
then $c\in\left(  \phi\left(  N\right)  :_{R}abL\right)  $. Now, assume that
$abcL\nsubseteq\phi\left(  N\right)  $. Thus by part (3) we have that
$cL\subseteq N$, as $abL\nsubseteq N$. Hence, $c\in\left(  N:_{R}L\right)  $.
Consequently, $\left(  N:_{R}abL\right)  =\left(  \phi\left(  N\right)
:_{R}abL\right)  \cup\left(  N:_{R}L\right)  $ and then $\left(
N:_{R}abL\right)  =\left(  \phi\left(  N\right)  :_{R}abL\right)  $ or
$\left(  N:_{R}abL\right)  =\left(  N:_{R}L\right)  $.

$(4)\Longrightarrow(5)$ Let for some nonunits $a,b\in R$, an ideal $I$ of $R$
and submodule $L$ of $M$ with $abIL\subseteq N$ and $abIL\nsubseteq\phi\left(
N\right)  $. Hence $I\subseteq\left(  N:_{R}abL\right)  $ and $I\nsubseteq$
$\left(  \phi\left(  N\right)  :_{R}abL\right)  $. If \ $abL\subseteq N$, then
we are done. So, assume that $abL\nsubseteq N$. Therefore by part (4) we have
$I\subseteq\left(  N:_{R}abL\right)  =\left(  N:_{R}L\right)  $. Thus,
$IL\subseteq N$.

$(5)\Longrightarrow(6)$ Suppose that $aIJL\subseteq N$, $aIJL\nsubseteq
\phi\left(  N\right)  $ and $JL\nsubseteq N$. Now, take $x\in I$. Since
$aIJL\subseteq N$ and $aIJL\nsubseteq\phi\left(  N\right)  $, there exists
$y\in I$ such that $ayJL\subseteq N$ $ayJL\nsubseteq\phi\left(  N\right)  $
and then by part (5), $ayL\subseteq N$. If $axJL\subseteq N$ $axJL\nsubseteq
\phi\left(  N\right)  $ and then similarly we have $axL\subseteq N$. Then we
have $a\left(  x+y\right)  JL\subseteq N$ $a\left(  x+y\right)  JL\nsubseteq
\phi\left(  N\right)  $ giving that $a\left(  x+y\right)  L\subseteq N$ which
implies that $axL\subseteq N$. Hence, we conclude that $aIL\subseteq N$.

$(6)\Longrightarrow(7)$ Suppose that $IJL\nsubseteq N$ and $a\in\left(
N:_{R}IJL\right)  $. Then $aIJL\subseteq N$ and there exists $x\in J$ such
that $xIL\nsubseteq N$. If $aIJL\subseteq\phi\left(  N\right)  $, then
$a\in\left(  \phi\left(  N\right)  :_{R}IJL\right)  $. Now, assume that
$aIJL\nsubseteq\phi\left(  N\right)  $. Then there exists $y\in J$ such that
$ayIL=yI\left(  Ra\right)  L\subseteq N$ and $ayIL=yI\left(  Ra\right)
\subseteq\phi\left(  N\right)  $. Then by part (6) we conclude that
$yIL\subseteq N$ or $\left(  Ra\right)  L\subseteq\phi\left(  N\right)  $. In
the latter case, we have $a\in\left(  N:_{R}L\right)  $. Assume tha former
case $yIL\subseteq N$. If $axIL=xI\left(  Ra\right)  L\subseteq N$ and
$axIL=xI\left(  Ra\right)  L\nsubseteq\phi\left(  N\right)  $, then by part
(6), we have $\left(  Ra\right)  \subseteq N$ which implies that $a\in\left(
N:_{R}L\right)  $. So we may assume that $axIL\subseteq\phi\left(  N\right)
$. We conclude that $x+y$ is a nonunit, since $J$ is a proper ideal. In this
case, we may assume that $a\left(  x+y\right)  IL=\left(  x+y\right)  I\left(
Ra\right)  L\subseteq N$ and $a\left(  x+y\right)  IL=\left(  x+y\right)
I\left(  Ra\right)  L\nsubseteq\phi\left(  N\right)  $. Again by part (6), we
conclude that $\left(  x+y\right)  IL\subseteq N$ or $\left(  Ra\right)
L\subseteq N$. If $\left(  x+y\right)  IL\subseteq N$, then we get
$xIL\subseteq N$, as $yIL\subseteq N$. This is a contradiction. Thus we
conclude that $\left(  Ra\right)  L\subseteq N$, that is, $a\in\left(
N:_{R}L\right)  $. By above arguments, we conclude that $\left(
N:_{R}IJL\right)  \subseteq\left(  \phi\left(  N\right)  :_{R}IJL\right)
\cup\left(  N:_{R}L\right)  $. Since reverse inclusion is always true, we have
the equality $\left(  N:_{R}IJL\right)  =\left(  \phi\left(  N\right)
:_{R}IJL\right)  \cup\left(  N:_{R}L\right)  $. In this case, we have $\left(
N:_{R}IJL\right)  =\left(  \phi\left(  N\right)  :_{R}IJL\right)  $ or
$\left(  N:_{R}IJL\right)  =\left(  N:_{R}L\right)  $.

$(7)\Longrightarrow(8)$ Suppose that $IJKL\subseteq N$ and $IJKL\nsubseteq
\phi\left(  N\right)  $. Then we have $K\subseteq\left(  N:_{R}IJL\right)  $
and $K\nsubseteq\left(  \phi\left(  N\right)  :_{R}IJL\right)  $. Then by part
(7), we have $K\subseteq\left(  N:_{R}IJL\right)  =(N:_{R}L)$ which implies
that $KL\subseteq N$.

$(8)\Longrightarrow(1)$ Is trivial.
\end{proof}

\begin{theorem}
\label{theo9}$M$ be a multiplication $R$-module and $N$ be a proper submodule
of $M$. Suppose that $\phi:S(M)\longrightarrow S(M)\cup\left\{  \emptyset
\right\}  $ is a function. Then the following conditions are equivalent:

\begin{enumerate}
\item $N$ is a $\phi-$classical 1-absorbing prime submodule of $M$;

\item If $N_{1}N_{2}N_{3}N_{4}\subseteq N$ for some proper submodules
$N_{1},N_{2},N_{3}$ and for some submodule $N_{4}$ of $M$ \ such that
$N_{1}N_{2}N_{3}N_{4}\nsubseteq\phi\left(  N\right)  $, then either
$N_{1}N_{2}N_{4}\subseteq N$ or $N_{3}N_{4}\subseteq N$.
\end{enumerate}
\end{theorem}

\begin{proof}
$(1)\Longrightarrow(2)$ Let $\ N_{1}N_{2}N_{3}N_{4}\subseteq N$ for some
proper submodules $N_{1},N_{2},N_{3}$ $,N_{4}$ of $M$ such that $N_{1}%
N_{2}N_{3}N_{4}\nsubseteq\phi\left(  N\right)  $. Since $M$ is multiplication
$R$-module,there are proper ideals $I_{1},I_{2},I_{3}$ of $R$ such that
$N_{1}=I_{1}M$ , $N_{2}=I_{2}M$ , $N_{3}=I_{3}M$. Therefore $I_{1}I_{2}%
I_{3}N_{4}\subseteq N$ and $I_{1}I_{2}I_{3}N_{4}\nsubseteq\phi\left(
N\right)  $, so by Theorem \ref{theo7}, we have $I_{1}I_{2}N_{4}\subseteq N$
or $I_{3}N_{4}\subseteq N$. Hence, $N_{1}N_{2}N_{4}\subseteq N$ or $N_{3}%
N_{4}\subseteq N$.

$(2)\Longrightarrow(1)$ Put $N_{1}=aM$, $N_{2}=bM$, $N_{3}=cM$ where $a,b,c$
are nonunit elements of $R$ and $N_{4}=Rm$, then the claim is clear.
\end{proof}

Recall from \cite{X14} that a $R-$module $F$ is said to be a flat $R-$module
if for each exact sequence $K\longrightarrow L\longrightarrow M$ of
$R-$modules, the sequence $F\otimes K\longrightarrow F\otimes L\longrightarrow
F\otimes M$ is also exact. Also, $F$ is said to be a faithfully flat the
sequence $K\longrightarrow L\longrightarrow M$ is exact if and only if
$F\otimes K\longrightarrow F\otimes L\longrightarrow F\otimes M$ is exact.
Azizi in \cite[lemma 3.2]{X4} showed that if $M$ is a $R-$module, $N$ is a
submodule of $M$ and $F$ is a flat $R-$module, then $(F\otimes N:_{F\otimes
M}a)=F\otimes(N:_{M}a)$ for every $a\in R$.

\begin{theorem}
\label{theo10}Let $R$ be a um-ring and $M$ \ be an $R$-module and
$\phi:S(M)\longrightarrow S(M)\cup\left\{  \emptyset\right\}  $ be a function.
Suppose that $N$ is a proper submodule of $M$ such that $F\otimes\phi\left(
N\right)  =\phi\left(  F\otimes N\right)  $.

\begin{enumerate}
\item If $F$ is flat $R$-module and $N$ is a $\phi-$ classical 1-absorbing
prime submodule of $M$ such that $F\otimes N\neq F\otimes M$, then $F\otimes
N$ is a $\phi-$ classical 1-absorbing prime submodule of $F\otimes M$.

\item Suppose that $F$ is a faithfully flat $R$-module. Then $N$ is a $\phi-$
classical 1-absorbing prime submodule of $M$ if and only if $F\otimes N$ is a
$\phi-$ classical 1-absorbing prime submodule of $F\otimes M$.
\end{enumerate}
\end{theorem}

\begin{proof}
$(1)$ Let $a,b,c\in R$ be nonunits. Then by Theorem \ref{theo7}, either
$\left(  N:_{M}abc\right)  =\left(  \phi\left(  N\right)  :_{M}abc\right)  $
or $\left(  N:_{M}abc\right)  =\left(  N:_{M}ab\right)  $ or $\left(
N:_{M}abc\right)  =\left(  N:_{M}c\right)  $. Assume that $\left(
N:_{M}abc\right)  =\left(  \phi\left(  N\right)  :_{M}abc\right)  $ . Then by
\cite[lemma 3.2]{X4} , $\left(  F\otimes N:_{F\otimes M}abc\right)
=F\otimes\left(  N:_{M}abc\right)  =F\otimes\left(  \phi\left(  N\right)
:_{M}abc\right)  =\left(  F\otimes\phi\left(  N\right)  :_{F\otimes
M}abc\right)  =\left(  \phi\left(  F\otimes N\right)  :_{F\otimes
M}abc\right)  $. Now, suppose that $\left(  N:_{M}abc\right)  =\left(
N:_{M}ab\right)  $. Again by \cite[lemma 3.2]{X5}, $\left(  F\otimes
N:_{F\otimes M}abc\right)  =F\otimes\left(  N:_{M}abc\right)  =F\otimes\left(
N:_{M}ab\right)  =\left(  F\otimes N:_{F\otimes M}ab\right)  $. Similarly, we
can show that if\ $\left(  N:_{M}abc\right)  =\left(  N:_{M}c\right)  $, then
$\left(  F\otimes N:_{F\otimes M}abc\right)  =$ $F\otimes\left(
N:_{M}abc\right)  =F\otimes\left(  N:_{M}c\right)  =\left(  F\otimes
N:_{F\otimes M}c\right)  $. Consequently, by Theorem \ref{theo7} we deduce
that $F\otimes N$ is a $\phi-$ classical 1-absorbing prime submodule of
$F\otimes M$.

$(2)$ Let $N$ be a $\phi-$ classical 1-absorbing prime submodule of $M$
$\ $\ and assume $F\otimes N=F\otimes M$. Then $0\rightarrow F\otimes
N\overset{\subseteq}{\longrightarrow}F\otimes M\rightarrow0$ is an exact
sequence. Since $F$ is a faithfully flat $R$-module, $0\rightarrow
N\overset{\subseteq}{\longrightarrow}M\rightarrow0$ is an exact sequence .
This implies that $N=M$ , which is a contradiction. So $F\otimes N\neq
F\otimes M$. Then $F\otimes N$ is a $\phi-$ classical 1-absorbing prime
submodule of $F\otimes M$ by part(1). Now for the converse , let $F\otimes N$
be a $\phi-$classical 1-absorbing prime submodule of $F\otimes M$. We have
$F\otimes N\neq F\otimes M$ and so $N\neq M$. Let $a,b,c\in R$ be nonunits.
Then by Theorem \ref{theo7} $\left(  F\otimes N:_{F\otimes M}abc\right)
=\left(  \phi\left(  F\otimes N\right)  :_{F\otimes M}abc\right)  $ or
$\left(  F\otimes N:_{F\otimes M}abc\right)  =\left(  F\otimes N:_{F\otimes
M}ab\right)  $ or $\left(  F\otimes N:_{F\otimes M}abc\right)  =\left(
F\otimes N:_{F\otimes M}c\right)  $ . Suppose that $\left(  F\otimes
N:_{F\otimes M}abc\right)  =\left(  \phi\left(  F\otimes N\right)  :_{F\otimes
M}abc\right)  $. Hence $F\otimes\left(  N:_{M}abc\right)  =\left(  F\otimes
N:_{F\otimes M}abc\right)  =\left(  \phi\left(  F\otimes N\right)  :_{F\otimes
M}abc\right)  =\left(  F\otimes\phi\left(  N\right)  :_{F\otimes M}abc\right)
=F\otimes\left(  \phi\left(  N\right)  :_{M}abc\right)  $. Thus $0\rightarrow
F\otimes\left(  \phi\left(  N\right)  :_{M}abc\right)  \overset{\subseteq
}{\longrightarrow}F\otimes\left(  N:_{M}abc\right)  \rightarrow0$ is an exact
sequence. Since $F$ is a faithfully flat $R$-module $0\rightarrow\left(
\phi\left(  N\right)  :_{M}abc\right)  \overset{\subseteq}{\longrightarrow
}\left(  N:_{M}abc\right)  \rightarrow0$ is an exact sequence which implies
that $\left(  N:_{M}abc\right)  =\left(  \phi\left(  N\right)  :_{M}%
abc\right)  $. Now, assume that $\left(  F\otimes N:_{F\otimes M}abc\right)
=\left(  F\otimes N:_{F\otimes M}ab\right)  $. Hence $F\otimes\left(
N:_{M}abc\right)  =\left(  F\otimes N:_{F\otimes M}abc\right)  =\left(
F\otimes N:_{F\otimes M}ab\right)  =F\otimes\left(  N:_{M}ab\right)  $. So
$0\rightarrow F\otimes\left(  N:_{M}ab\right)  \overset{\subseteq
}{\longrightarrow}F\otimes\left(  N:_{M}abc\right)  \rightarrow0$ is an exact
sequence. Since $F$ is a faithfully flat $R$-module, $0\rightarrow\left(
N:_{M}ab\right)  \overset{\subseteq}{\longrightarrow}\left(  N:_{M}abc\right)
\rightarrow0$ is an exact sequence which implies that $\left(  N:_{M}%
abc\right)  =\left(  N:_{M}ab\right)  $. In other case, one can similarly show
that $\left(  N:_{M}abc\right)  =\left(  N:_{M}c\right)  $.Suppose that
$\left(  F\otimes N:_{F\otimes M}abc\right)  =\left(  F\otimes N:_{F\otimes
M}c\right)  $. Hence $F\otimes\left(  N:_{M}abc\right)  =\left(  F\otimes
N:_{F\otimes M}abc\right)  =\left(  F\otimes N:_{F\otimes M}c\right)
=F\otimes\left(  N:_{M}c\right)  $. Hence, $0\rightarrow F\otimes\left(
N:_{M}c\right)  \overset{\subseteq}{\longrightarrow}F\otimes\left(
N:_{M}abc\right)  \rightarrow0$ is an exact sequence. Again since $F$ is a
faithfully flat $R$-module, $0\rightarrow\left(  N:_{M}c\right)
\overset{\subseteq}{\longrightarrow}\left(  N:_{M}abc\right)  \rightarrow0$ is
an exact sequence which implies that $\left(  N:_{M}abc\right)  =\left(
N:_{M}c\right)  $. Consequently, $N$ is a $\phi-$ classical 1-absorbing prime
submodule of $M$ by Theorem \ref{theo7}.
\end{proof}

\begin{corollary}
\label{cor4}Let $R$ be a $um$-ring, $M$ \ be an $R$-module and $X$ be an
indeterminate. If $N$ is a $\phi-$ classical 1-absorbing prime submodule of
$M$ with $R\left[  X\right]  \otimes N=\phi\left(  R\left[  X\right]  \otimes
N\right)  $, then $N\left[  X\right]  $ is a $\phi-$ classical 1-absorbing
prime submodule of $M\left[  X\right]  $.
\end{corollary}

\begin{proof}
Assume that $N$ is a $\phi-$classical 1-absorbing prime submodule of $M$ with
$R\left[  X\right]  \otimes N=\phi\left(  R\left[  X\right]  \otimes N\right)
$. Notice that $R\left[  X\right]  $ is a flat $R$-module. So by Theorem
\ref{theo10}, $R\left[  X\right]  \otimes N$ is a $\phi-$ classical
1-absorbing prime submodule of $R\left[  X\right]  \otimes M$. Since,
$R\left[  X\right]  \otimes N\cong N\left[  X\right]  $, $R\left[  X\right]
\otimes M\cong M\left[  X\right]  $ and then $N\left[  X\right]  $ is a
$\phi-$classical \ 1-absorbing prime submodule of $M\left[  X\right]  $.
\end{proof}

\begin{theorem}
\label{theo12}Let $M$ be an $R-$module and a be an element of $R$ such that
$aM\neq M$. Suppose that $\left(  0:_{M}a\right)  \subseteq aM$. Then $aM$ is
an almost classical 1-absorbing prime submodule of $M$ if and only if it is a
classical 1-absorbing prime submodule of $M$.
\end{theorem}

\begin{proof}
Assume that $aM$ is an almost classical 1-absorbing prime submodule of $M$.
Let nonunits $x,y,z\in R$ and $m\in M$ such that $\ xyzm\in aM$. We will show
that $xym\in aM$ or $zm\in aM$. If $xyzm\notin(aM:_{R}M)aM$, then there is
nothing to prove. So, assume that $xyzm\in(aM:_{R}M)aM$. Note that
$\ xy(z+a)m\in aM$. If $z+a$ is unit, then we have $xym\in aM$, as required.
So, assume that $z+a$ is a nonunit. If $\ xy(z+a)m\notin(aM:_{R}M)aM$, then
$xym\in aM$ or $\left(  z+a\right)  m\in aM$ and we are done. Therefore,
suppose that $\ xy(z+a)m\in(aM:_{R}M)aM$. Hence $\ xyzm\in(aM:_{R}M)aM$ gives
$\ xyam\in(aM:_{R}M)aM$. Then, there exists $m%
\acute{}%
\in(aM:_{R}M)M$ such that $xyam=am%
\acute{}%
$ and so $xym-m%
\acute{}%
\in\left(  0:_{M}a\right)  \subseteq aM$ which shows $xym\in aM$, as $m%
\acute{}%
\in aM$. Consequently, $aM$ is a classical 1-absorbing prime. The converse is
easy to check.
\end{proof}

Let $S$ be a multiplicatively closed subset of $R$. It is well-known that each
submodule of $S^{-1}M$ is in the form of $S^{-1}N$ for some submodule $N$ of
$M$. Let $\phi:S(M)\longrightarrow S(M)\cup\left\{  \emptyset\right\}  $ be a
function and define $\phi_{S}:S(S^{-1}M)\longrightarrow S(S^{-1}M)\cup\left\{
\emptyset\right\}  $ by $\phi_{S}(S^{-1}N)=S^{-1}\phi\left(  N\right)  $ (and
$\phi_{S}(S^{-1}N)=\emptyset$ when $\phi\left(  N\right)  =\emptyset$) for
every submodule $N$ of $M$. For an $R-$module $M,$ the set of zero-divisors of
$M$ is denoted by $Z_{R}\left(  M\right)  $.

\begin{theorem}
\label{theo13}Let $M$ be an $R-$module, $S$ be a multiplicatively closed
subset of $R$ such that $(N:_{R}M)\cap S=\emptyset$ and $N$ be a $\phi
-$classical 1-absorbing prime submodule of $M$. Then $S^{-1}N$ is a $\phi
_{S}-$classical 1-absorbing prime submodule of $S^{-1}M$.
\end{theorem}

\begin{proof}
Let $N$ be a $\phi-$ classical 1-absorbing prime submodule of $M$ and $\left(
N:_{R}M\right)  \cap S=\emptyset$. Suppose that $\frac{a_{1}}{s_{1}}%
\frac{a_{2}}{s_{2}}\frac{a_{3}}{s_{3}}\frac{m}{s_{4}}\in S^{-1}N\setminus$
$\phi_{S}(S^{-1}N)$ for some nonunits $\frac{a_{1}}{s_{1}},\frac{a_{2}}{s_{2}%
},\frac{a_{3}}{s_{3}}\in S^{-1}R$ and $\frac{m}{s_{4}}\in S^{-1}M$. Then there
exists $s^{\ast}\in S$ such that $s^{\ast}(a_{1}a_{2}a_{3}m)=a_{1}a_{2}%
a_{3}(sm)\in N$. If $sa_{1}a_{2}a_{3}m\in\phi\left(  N\right)  $, then
$\frac{a_{1}}{s_{1}}\frac{a_{2}}{s_{2}}\frac{a_{3}}{s_{3}}\frac{m}{s_{4}%
}=\frac{sa_{1}a_{2}a_{3}m}{ss_{1}s_{2}s_{3}s_{4}}\in S^{-1}\phi\left(
N\right)  =\phi_{S}(S^{-1}N)$, a contradiction. Since $N$ is a $\phi
-$classical prime submodule, then we have $a_{1}a_{2}(sm)\in N$ or
$a_{3}\left(  sm\right)  \in N$.Thus $\frac{a_{1}}{s_{1}}\frac{a_{2}}{s_{2}%
}\frac{m}{s_{4}}=\frac{sa_{1}a_{2}m}{ss_{1}s_{2}s_{4}}\in S^{-1}N$ or
$\frac{a_{3}}{s_{3}}\frac{m}{s_{4}}=\frac{sa_{3}m}{ss_{3}s_{4}}\in S^{-1}N$.
Consequently, $S^{-1}N$ is a $\phi_{S}-$classical 1-absorbing prime submodule
of $S^{-1}M$.
\end{proof}

\begin{definition}
Let $N$ be a proper submodule of $M$ and $a,b,c$ be nonunits in $R$ and $m\in
M$. If $N$ is a $\phi-$ classical 1-absorbing prime submodule and $abcm\in
\phi\left(  N\right)  $, $abm\notin N$ and $\ cm\notin N$, then we say that
$\left(  a,b,c,m\right)  $ is called a $\phi-$classical 1-quadruple-zero of
$N$.
\end{definition}

\begin{theorem}
\label{theo11}Let $N$ be a $\phi-$classical 1-absorbing prime submodule of an
$R-$module $M$ and suppose that $abcK\subseteq N$ for some nonunits $a,b,c\in
R$ and some submodule $K$ of $M$. If $\left(  a,b,c,k\right)  $ is not a
$\phi-$classical 1-quadruple-zero of $N$ for any $k\in K$, then $abK\subseteq
N$ or $\ cK\subseteq N$.
\end{theorem}

\begin{proof}
Suppose that $\left(  a,b,c,k\right)  $ is not a $\phi-$classical
1-quadruple-zero of $N$ for every $k\in K$. Assume on the contrary that
$abK\nsubseteq N$ and $cK\nsubseteq N$. Then there are $k_{1},k_{2}\in K$ such
that $abk_{1}\notin N$ and $ck_{2}\notin N$. If $abck_{1}\notin\phi\left(
N\right)  $, then we have $ck_{1}\in N$. If $abck_{1}\in\phi\left(  N\right)
$, then since $abk_{1}\notin N$ and $\left(  a,b,c,k_{1}\right)  $ is not a
$\phi-$ classical 1-quadruple-zero of $N$, we conclude again that $ck_{1}\in
N$. By a similar argument, since $\left(  a,b,c,k_{2}\right)  $ is not a
$\phi-$classical 1-quadruple-zero of $N$ and $ck_{2}\notin N$, then we deduce
that $abk_{2}\in N$. By our hypothesis, $abc(k_{1}+k_{2})\in N$ and $\left(
a,b,c,k_{1}+k_{2}\right)  $ is not a $\phi-$classical 1-quadruple-zero of $N$.
Hence we have either $ab(k_{1}+k_{2})\in N$ or $c(k_{1}+k_{2})\in N.$ If
$ab(k_{1}+k_{2})\in N=abk_{1}+abk_{2}\in N$, then since $abk_{2}\in N$, we
have $abk_{1}\in N$, which is a contradiction. If $c(k_{1}+k_{2}%
)=ck_{1}+ck_{2}\in N$, then since $ck_{1}\in N,$ we have $ck_{2}\in N$, which
again is a contradiction. Thus, we deduce $abK\subseteq N$ or $\ cK\subseteq
N$.
\end{proof}

We introduce the following definition to give a further characterization of
this class of submodules in terms of some ideals and submodules.

\begin{definition}
Let $N$ be a $\phi-$classical 1-absorbing prime submodule of an $R$-module $M$
and suppose that $HIJK\subseteq N$ for some ideals $H,I,J$ of $R$ and some
submodule $K$ of $M$. We say that $N$ is a free $\phi-$classical
1-quadruple-zero with respect to $HIJK$ if $\left(  a,b,c,k\right)  $ is not a
$\phi-$classical 1-quadruple-zero of $N$ for every $a\in H$, $b\in I$, $c\in
J$ and $k\in K$ .
\end{definition}

Let $N$ be a $\phi-$classical 1-absorbing prime submodule of $M$ and suppose
that $HIJK\subseteq N$ for some ideals $H,I,J$ of $R$ and some submodule $K$
of $M$ such that $N$ is a free $\phi-$classical 1-quadruple-zero with respect
to $HIJK$. Hence, if $a\in H$, $b\in I$, $c\in J$ and $k\in K$, then $abk\in
N$ or $ck\in N$.

\begin{corollary}
\label{cor5}Let $N$ be a $\phi-$ classical 1-absorbing prime submodule of an
$R$-module $M$ and suppose that $HIJK\subseteq N$ for some ideals $H,I,J$ of
$R$ and some submodule $K$ of $M$. If $N$ is a free $\phi-$classical
1-quadruple-zero with respect to $HIJK$, then $HIK\subseteq N$ or $JK\subseteq
N$.
\end{corollary}

\begin{proof}
Suppose that $N$ is a free $\phi-$classical 1-quadruple-zero with respect to
$HIJK$. Assume that $HIK\nsubseteq N$ and $JK\nsubseteq N$. Then there are
$a\in H$, $b\in I$, $c\in J$ with $abK\nsubseteq N$ and $cK\nsubseteq N$.
Since $abcK\subseteq N$ and $N$ is a free $\phi-$classical 1-quadruple-zero
with respect to $HIJK$, then Theorem \ref{theo11} implies that $abK\subseteq
N$ and $cK\subseteq N$, which is a contradiction. Therefore, $HIK\subseteq N$
or $JK\subseteq N$.
\end{proof}

\begin{theorem}
\label{theo14}Let $N$ be a $\phi-$classical 1-absorbing prime submodule of $M$
and suppose that $\left(  a,b,c,m\right)  $ is a $\phi-$classical
1-quadruple-zero of $N$ for some nonunits $a,b,c\in R$ and $m\in M$. Then

\begin{enumerate}
\item $abcN\subseteq\phi\left(  N\right)  $.

\item $ab(N:_{R}M)m\subseteq\phi\left(  N\right)  $.

\item If $a,b\notin(N:_{R}cm)$, then $ac(N:_{R}M)m\subseteq\phi\left(
N\right)  $ and $bc(N:_{R}M)\subseteq\phi\left(  N\right)  $.

\item If $a,b\notin(N:_{R}cm)$, then $a(N:_{R}M)^{2}m\subseteq\phi\left(
N\right)  $, $b(N:_{R}M)^{2}m\subseteq\phi\left(  N\right)  $ and
$c(N:_{R}M)^{2}m\subseteq\phi\left(  N\right)  $.

\item If $a,b\notin(N:_{R}cm)$, then $(N:_{R}M)^{3}m\subseteq\phi\left(
N\right)  $.
\end{enumerate}
\end{theorem}

\begin{proof}
$(1)$ Suppose that $abcN\nsubseteq\phi\left(  N\right)  $. Then there exists
$n\in N$ with $abcn\notin\phi\left(  N\right)  $. Hence $abc(m+n)\in
N\setminus\phi\left(  N\right)  $, so we conclude that $ab(m+n)\in N$ or
$c(m+n)\in N$. Thus $abm\in N$ or $\ cm\in N$, which contradicts the
assumption that $\left(  a,b,c,m\right)  $ is $\phi-$classical
1-quadruple-zero. Thus $abcN\subseteq\phi\left(  N\right)  $.

$(2)$ Assume that $abxm\notin\phi\left(  N\right)  $ for some $x\in(N:_{R}M)$.
Then $ab(c+x)m\notin\phi\left(  N\right)  $. Since $xm\in N$, $ab(c+x)m\in N$.
If $c+x$ is a unit, then $abm\in N$ which is a contradiction. Thus $c+x$ is
nonunit and since $N$ is a $\phi-$classical 1-absorbing prime submodule we
have either $abm\in N$ or $\left(  c+x\right)  m\in N$. Thus $abm\in N$ or
$cm\in N$, again a contradiction. Thus, $ab(N:_{R}M)m\subseteq\phi\left(
N\right)  $.

$(3)$ Let $a,b\notin(N:_{R}cm)$, that is $acm\notin N$ and $bcm\notin N$.
Choose $x\in(N:_{R}M)$ and $acxm\notin\phi\left(  N\right)  $. Then we have
$a\left(  b+x\right)  cm\in N\setminus\phi\left(  N\right)  $. Note that $b+x$
must be nonunit. As $N$ is $\phi-$classical 1-absorbing prime, we have
$a\left(  b+x\right)  m\in N$ or $cm\in N$ which yields that $abm\in N$ or
$cm\in N$. This is a contradiction. Hence, we have $ac(N:_{R}M)m\subseteq
\phi\left(  N\right)  $. Likewise, we have $bc(N:_{R}M)m\subseteq\phi\left(
N\right)  $.

$(4)$ Let $a,b\notin(N:_{R}cm)$, that is $acm\notin N$ and $bcm\notin N$. We
want to show $a(N:_{R}M)^{2}m\subseteq\phi\left(  N\right)  $ and to see this,
suppose the contrary. Then there exist \ $x,y\in(N:_{R}M)$ such that
$axym\notin\phi\left(  N\right)  $. Since $x,y\in(N:_{R}M)$, then $a\left(
b+x\right)  \left(  c+y\right)  m\in N$ and $a\left(  b+x\right)  \left(
c+y\right)  m\notin\phi\left(  N\right)  $, as $abcm\in\phi\left(  N\right)
$. Also, note that $b+x$ and $c+y$ are nonunits, as $acm\notin N$ and
$abm\notin N$. Since $N$ is a $\phi-$ classical 1-absorbing prime submodule,
we have either $a\left(  b+x\right)  m\in N$ or $\left(  c+y\right)  m\in N$.
Which implies that $abm\in N$ or $cm\in N$. This is a contradiction. Thus, we
have $a(N:_{R}M)^{2}m\subseteq\phi\left(  N\right)  $. Similarly we can prove
that $b(N:_{R}M)^{2}m\subseteq\phi\left(  N\right)  $ and $c(N:_{R}%
M)^{2}m\subseteq\phi\left(  N\right)  $.

$\left(  5\right)  $ Let $a,b\notin(N:_{R}cm)$, that is $acm\notin N$ and
$bcm\notin N$. Suppose that $(N:_{R}M)^{3}m\nsubseteq\phi\left(  N\right)  $.
Then there exist $x,y,z\in(N:_{R}M)$ such that $xyzm\notin\phi\left(
N\right)  $. Then we have $\left(  a+x\right)  \left(  b+y\right)  \left(
c+z\right)  m\in N\setminus\phi\left(  N\right)  $ from by parts (2), (3) ve
(4). Also, since $bcm\notin N$ and $abm\notin N$ we get $\left(  a+x\right)
$, $\left(  b+y\right)  $ and $\left(  c+z\right)  $ are nonunits. Then we
conclude that $\left(  a+x\right)  \left(  b+y\right)  m\in N$ or $\left(
c+z\right)  m\in N$. This implies that $abm\in N$ or $cm\in N$, which is a
contradiction. Hence, $(N:_{R}M)^{3}m\subseteq\phi\left(  N\right)  $.
\end{proof}

\begin{theorem}
\label{theo15}If $N$ is a $\phi-$classical 1-absorbing prime submodule of an
$R$-module $M$ that is not classical 1-absorbing prime, then there exists
$\phi-$classical 1-quadruple-zero $\left(  a,b,c,m\right)  $ of $N$. If
$a,b\notin(N:_{R}M)$ and then $(N:_{R}M)^{3}N\subseteq\phi\left(  N\right)  $.
\end{theorem}

\begin{proof}
Assume that $(N:_{R}M)^{3}N\nsubseteq\phi\left(  N\right)  $, then there are
$x_{1,}x_{2,}x_{3}\in(N:_{R}M)$ and $n\in N$ such that $x_{1}x_{2}x_{3}%
n\notin\phi\left(  N\right)  $. By Theorem \ref{theo14}, $(a+x_{1}%
)(b+x_{2})(c+x_{3})(m+n)\in N\setminus\phi\left(  N\right)  $. Since $N$ is a
$\phi-$classical 1-absorbing prime submodule we have that either
$(a+x_{1})(b+x_{2})(m+n)\in N$ or $(c+x_{3})(m+n)\in N$. Therefore, $abm\in N$
or $cm\in N$, a contradiction. This shows that $(N:_{R}M)^{3}N\subseteq
\phi\left(  N\right)  $.
\end{proof}

\begin{corollary}
\label{cor6}Let $M$ be an $R-$module and $N$ be a $\phi-$classical 1-absorbing
prime submodule of $M$ such that $\phi\left(  N\right)  \subseteq(N:_{R}%
M)^{4}N$. Then $N$ is $\omega-$classical 1-absorbing prime.
\end{corollary}

\begin{proof}
If \ $N$ is a classical 1-absorbing prime submodule of $M$, then it is clear.
Hence, suppose that $N$ is not a classical 1-absorbing prime submodule of $M$.
Therefore by Theorem \ref{theo15} we have $(N:_{R}M)^{3}N\subseteq\phi\left(
N\right)  \subseteq(N:_{R}M)^{4}N\subseteq(N:_{R}M)^{3}N$, that is,
$\phi\left(  N\right)  =(N:_{R}M)^{3}N=(N:_{R}M)^{4}N$. Therefore,
$\phi\left(  N\right)  =(N:_{R}M)^{j}N$ for all $j\geq3$ and the result is obtained.
\end{proof}

As a direct consequence of Theorem \ref{theo15} we have the following result.

\begin{corollary}
\label{cor7}Let $M$ be an $R-$module and $N$ be a proper submodule of $M$. If
$N$ is an n-almost classical 1-absorbing prime submodule $\left(
n\geq4\right)  $ of $M$ that is not classical 1-absorbing prime, then
$(N:_{R}M)^{3}N=(N:_{R}M)^{n-1}N$.
\end{corollary}

\begin{corollary}
\label{cor8}Let $M$ be a multipllication $R-$module and $N$ be a proper
submodule of $M$.

\begin{enumerate}
\item If $N$ is a $\phi-$classical 1-absorbing prime submodule of $M$ that is
not classical 1-absorbing prime, then $N^{4}\subseteq\phi\left(  N\right)  $.

\item If $N$ is an n-almost classical 1-absorbing prime submodule $\left(
n\geq4\right)  $ of $M$ that is not classical 1-absorbing prime, then
$N^{4}\subseteq N^{n}$.
\end{enumerate}
\end{corollary}

\begin{proof}
$\left(  1\right)  $Since $M$ is multiplication, then $N=(N:_{R}M)M$.
Therefore by Theorem \ref{theo15} and Remark \ref{remark1}, $N^{4}%
=(N:_{R}M)^{3}N\subseteq\phi\left(  N\right)  $.

$\left(  2\right)  $Notice that $\phi_{n}\left(  N\right)  =(N:_{R}%
M)^{n-1}N=N^{n}$. Now, use part (1).
\end{proof}

\begin{theorem}
\label{theo16}Let $N$ be a $\phi-$classical 1-absorbing prime submodule of
$M$. If $N$ is not classical 1-absorbing prime, then we have the following:

\begin{enumerate}
\item $\sqrt{(N:_{R}M)}=\sqrt{(\phi\left(  N\right)  :_{R}M)}$

\item If $M$ is multiplication, then $M-rad\left(  N\right)  =M-rad\left(
\phi\left(  N\right)  \right)  $.
\end{enumerate}
\end{theorem}

\begin{proof}
$\left(  1\right)  $Assume that $N$ is not classical 1-absorbing prime. By
Theorem \ref{theo15}, $(N:_{R}M)^{3}N\subseteq\phi\left(  N\right)  $. Then
$(N:_{R}M)^{4}=(N:_{R}M)^{3}(N:_{R}M)\subseteq\left(  (N:_{R}M)^{3}%
N:_{R}M\right)  \subseteq(\phi\left(  N\right)  :_{R}M)$ and so $(N:_{R}%
M)\subseteq(\phi\left(  N\right)  :_{R}M)$. Hence we have $\sqrt{(N:_{R}%
M)}=\sqrt{(\phi\left(  N\right)  :_{R}M)}$

$\left(  2\right)  $By part (1), $M-rad\left(  N\right)  =\sqrt{(N:_{R}%
M)}M=\sqrt{(\phi\left(  N\right)  :_{R}M)M=}M-rad\left(  \phi\left(  N\right)
\right)  $.
\end{proof}

\begin{theorem}
\label{theo17}Let $M$ be an $R-$module. Suppose that $N_{1},N_{2}$ are $\phi
-$classical 1-absorbing prime submodules of $M$ that are not classical
1-absorbing prime submodules. Then

\begin{enumerate}
\item $\sqrt{(N_{1}:_{R}M)+(N_{2}:_{R}M)}=\sqrt{(\phi\left(  N_{1}\right)
:_{R}M)+(\phi\left(  N_{2}\right)  :_{R}M)}$.

\item If $N_{1}+N_{2}\neq M$, $\phi\left(  N_{1}\right)  \subseteq N_{2}$ and
$\phi\left(  N_{2}\right)  \subseteq\phi\left(  N_{1}+N_{2}\right)  $, then
$N_{1}+N_{2}$ is a $\phi-$classical 1-absorbing prime submodule.
\end{enumerate}
\end{theorem}

\begin{proof}
$\left(  1\right)  $By Theorem \ref{theo16}, we have $\sqrt{(N_{1}:_{R}%
M)}=\sqrt[-]{(\phi\left(  N_{1}\right)  :_{R}M)}$ and $\sqrt{(N_{2}:_{R}%
M)}=\sqrt[-]{(\phi\left(  N_{2}\right)  :_{R}M)}$. Now, by \cite[$2.25\left(
i\right)  $]{X21} the result follows.

$\left(  2\right)  $Suppose that $N_{1}+N_{2}\neq M$, $\phi\left(
N_{1}\right)  \subseteq N_{2}$ and $\phi\left(  N_{2}\right)  \subseteq
\phi\left(  N_{1}+N_{2}\right)  $. Since $\left(  N_{1}+N_{2}\right)  \diagup
N_{2}\simeq N_{1}\diagup\left(  N_{1}\cap N_{2}\right)  $ and $N_{1}$ is
$\phi-$classical 1-absorbing prime, we get $\left(  N_{1}+N_{2}\right)
\diagup N_{2}$ is a weakly classical 1-absorbing prime submodule of $M\diagup
N_{2}$, by Theorem \ref{theo1}(3). Now, the assertion follows from Theorem
\ref{theo1}(4).
\end{proof}

Next, we will discuss the $\phi-$ classical 1-absorbing prime submodules of
cartesian product of modules.

\begin{theorem}
\label{theo18}Let $M_{1},M_{2}$ be $R$-modules and $N_{1}$ be a proper
submodule of $M_{1}$. Suppose that $\psi_{i}:S(M_{i})\longrightarrow
S(M_{i})\cup\left\{  \emptyset\right\}  $ be functions $\left(  \text{for
}i=1,2\right)  $ and let $\phi=\psi_{1}\times\psi_{2}$.Then the following
conditions are equivalent.

\begin{enumerate}
\item $N=N_{1}\times M_{2}$ is a $\phi-$ classical 1-absorbing prime submodule
of $M=M_{1}\times M_{2}$.

\item $N_{1}$ is a $\Psi_{1}-$ classical 1-absorbing prime submodule of
$M_{1}$ and for each nonunits $r,s,t\in R$ and $m_{1}\in M_{1}$ if
$rstm_{1}\in\psi_{1}\left(  N_{1}\right)  $, $rsm_{1}\notin N_{1}$,
$tm_{1}\notin N_{1}$, then $rst\in\left(  \psi_{2}(M_{2}):_{R}M_{2}\right)  $.
\end{enumerate}
\end{theorem}

\begin{proof}
$(1)\Longrightarrow(2)$ Suppose that $N=N_{1}\times M_{2}$ is a $\phi-$
classical 1-absorbing prime submodule of $M=M_{1}\times M_{2}$. Let nonunits
$r,s,t\in R$ and $m_{1}\in M_{1}$ be such that $rstm_{1}\in N_{1}\setminus
\psi_{1}\left(  N_{1}\right)  $. Then $rst\left(  m_{1},0\right)  \in
N\setminus\phi\left(  N\right)  $. Since $N$ is a $\phi-$classical 1-absorbing
prime submodule, we have either $rs\left(  m_{1},0\right)  \in N$ or $t\left(
m_{1},0\right)  \in N$, and so $rsm_{1}\in N_{1}$ or $tm_{1}\in N_{1}$.
Consequently, $N_{1}$ is a $\phi-$ classical 1-absorbing prime submodule of
$M_{1}$. Now, assume that $rstm_{1}\in\psi_{1}\left(  N_{1}\right)  $ for some
nonunits $r,s,t\in R$ and $m_{1}\in M_{1}$ such that $rsm_{1}\notin N_{1}$ and
$tm_{1}\notin N_{1}$. Suppose that $rst\notin\left(  \psi_{2}(M_{2}):_{R}%
M_{2}\right)  $. Therefore there exists $m_{2}\in M_{2}$ such that
$rstm_{2}\notin\psi_{2}(M_{2})$. Hence $rst\left(  m_{1},m_{2}\right)  \in
N\setminus\phi\left(  N\right)  $, and so $rs\left(  m_{1},m_{2}\right)  \in
N$ or $t\left(  m_{1},m_{2}\right)  \in N$ . Thus $rsm_{1}\in N_{1}$ or
$tm_{1}\in N_{1}$ which is a contradiction. Consequently, $rst\in\left(
\psi_{2}(M_{2}):_{R}M_{2}\right)  $.

$(2)\Longrightarrow(1)$ Let nonunits $r,s,t\in R$ and $\left(  m_{1}%
,m_{2}\right)  \in M=M_{1}\times M_{2}$ be such that $rst\left(  m_{1}%
,m_{2}\right)  \in N\setminus\phi\left(  N\right)  $. First assume that
$rstm_{1}\notin\psi_{1}\left(  N_{1}\right)  $. Then by part (2) $rsm_{1}\in
N_{1}$ or $tm_{1}\in N_{1}$. So, $rs\left(  m_{1},m_{2}\right)  \in N$ or
$t\left(  m_{1},m_{2}\right)  \in N$ and thus we are done. If $rstm_{1}\in
\psi_{1}\left(  N_{1}\right)  $, then $rstm_{2}\notin\psi_{2}(M_{2})$.
Therefore, $rst\notin\left(  \psi_{2}(M_{2}):_{R}M_{2}\right)  $ and so part
(2) implies that either $rsm_{1}\in N_{1}$ or $tm_{1}\in N_{1}$. Again we have
that $rs\left(  m_{1},m_{2}\right)  \in N$ or $t\left(  m_{1},m_{2}\right)
\in N$ which shows $N$ is a $\phi-$ classical 1-absorbing prime submodule of
$M$.
\end{proof}

\begin{corollary}
\label{cor9}Let $R=R_{1}\times R_{2}$ be a decomposable ring and
$M=M_{1}\times R_{2}$ be an $R-$module where $M_{1}$ is an $R_{1}$-module. If
$N_{1}$ is a weakly classical 1-absorbing prime submodule of $M_{1}$, then
$N=N_{1}\times R_{2}$ is a $4-$almost classical 1-absorbing prime submdule of
$M$.
\end{corollary}

\begin{proof}
Suppose that $N_{1}$ is a weakly classical 1-absorbing prime submodule of
$M_{1}$. If $N_{1}$ is a classical 1-absorbing prime submodule of $M_{1}$,
then it is easy to see that $N$ is a classical 1-absorbing prime submodule of
$M$ and so is $\phi-$classical 1-absorbing prime submodule of $M$, for all
$\phi$. Assume that $N_{1}$ is not classical 1-absorbing prime. Therefore by
Theorem \ref{theo15}, $(N_{1}:_{R_{1}}M_{1})^{3}N_{1}=\left\{  0\right\}  $
and so $\phi_{4}\left(  N\right)  =(N:_{R}M)^{3}N=\left\{  0\right\}  \times
R_{2}$.
\end{proof}

\begin{theorem}
\label{theo19}Let $R=R_{1}\times R_{2}$ be a decomposable ring and
$M=M_{1}\times M_{2}$ be an $R$-module where $M_{1}$ is an $R_{1}$-module and
$M_{2}$ is an $R_{2}$-module. Let $\psi_{i}:S(M_{i})\longrightarrow
S(M_{i})\cup\left\{  \emptyset\right\}  $ be functions for $i=1,2$ and let
$\phi=\psi_{1}\times\psi_{2}$. If $N=N_{1}\times M_{2}$ is a proper submodule
of $M$, then the following conditions are equivalent.

\begin{enumerate}
\item $N_{1}$ is a classical 1-absorbing prime submodule of $M_{1}$.

\item $N$ is a classical 1-absorbing prime submodule of $M$.

\item $N$ is $\phi-$ classical 1-absorbing prime submodule of $M$ where
$\psi_{2}\left(  M_{2}\right)  \neq M_{2\text{.}}$.
\end{enumerate}
\end{theorem}

\begin{proof}
$(1)\Longrightarrow(2)$ Let $(a_{1},a_{2})(b_{1},b_{2})(c_{1},c_{2}%
)(m_{1},m_{2})\in N$ for some nonunits $(a_{1},a_{2}),(b_{1},b_{2}%
),(c_{1},c_{2})\in R$ and $(m_{1},m_{2})\in M$. Then $\left(  a_{1}b_{1}%
c_{1}m_{1,}a_{2}b_{2}c_{2}m_{2}\right)  \in N=N_{1}\times M_{1}$. Hence
$a_{1}b_{1}c_{1}m_{1}\in N_{1}$. Since $N_{1}$ is a classical 1-absorbing
prime submodule of $M_{1}$we have $a_{1}b_{1}m_{1}\in N_{1}$ or $c_{1}m_{1}\in
N_{1}$ which shows either $(a_{1},a_{2})(b_{1},b_{2})(m_{1},m_{2})\in N$ or
$(c_{1},c_{2})(m_{1},m_{2})\in N$. Consequently, $N$ is a classical
1-absorbing prime submodule of $M$.

$(2)\Longrightarrow(3)$ It is clear that every classical 1-absorbing prime
submodule is a $\phi-$classical 1-absorbing prime submodule.

$(3)\Longrightarrow(1)$Let $abcm\in N_{1}$ for some nonunits $a,b,c\in R_{1}$
and $m_{1}\in M_{1}$. By assumption, there exists $m_{2}\in M_{2}\setminus$
$\psi_{2}\left(  M_{2}\right)  $. Thus $(a,1)(b,1)(c,1)(m,m_{2})\in
N\setminus\phi(N)$. So, we have $(a,1)(b,1)(m,m_{2})\in N$ or $(c,1)(m,m_{2}%
)\in N$. Hence, $\ abm\in N_{1}$ or $cm\in N_{1}$. Therefore, $N_{1}$ is a
classical 1-absorbing prime submodule of $M_{1}$.
\end{proof}

\begin{theorem}
\label{theo20}Let $R=R_{1}\times R_{2}$ be a decomposable ring and
$M=M_{1}\times M_{2}$ be an $R$-module where $M_{1}$ is an $R_{1}$-module and
$M_{2}$ is an $R_{2}$-module. $\psi_{i}:S(M_{i})\longrightarrow S(M_{i}%
)\cup\left\{  \emptyset\right\}  $ be functions for $i=1,2$ where $\psi
_{2}\left(  M_{2}\right)  =M_{2}$ and let $\phi=\psi_{1}\times\psi_{2}$. If
$N=N_{1}\times M_{2}$ is a proper submodule of $M$, then $N_{1}$ is a
$\psi_{1}-$classical 1-absorbing prime submodule of $M_{1}$ if and only if $N$
is $\phi-$classical 1-absorbing prime submodule of $M$.
\end{theorem}

\begin{proof}
Suppose that $N$ is a $\phi-$classical 1-absorbing prime submodule of $M$.
First we show that $N_{1}$ is a $\psi_{1}-$classical 1-absorbing prime
submodule of $M_{1}$ independently whether $\psi_{2}\left(  M_{2}\right)
=M_{2}$ or $\psi_{2}\left(  M_{2}\right)  \neq M_{2}$. Let $a_{1}b_{1}%
c_{1}m_{1}\in N_{1}\setminus\psi_{1}\left(  N_{1}\right)  $ for some nonunits
$a_{1},b_{1},c_{1}\in R_{1}$ and $m_{1}\in M_{1}$. Then $(a,1)(b,1)(c,1)(m_{1}%
,m_{2})\in\left(  N_{1}\times M_{2}\right)  \setminus\left(  \psi_{1}\left(
N_{1}\right)  \times\psi_{2}\left(  M_{2}\right)  \right)  =N\setminus
\phi\left(  N\right)  $ for any $m_{2}\in M_{2}$. Since $N$ is a $\phi
-$classical 1-absorbing prime submodule of $M$, we have either
$(a,1)(b,1)(m_{1},m_{2})\in\left(  N_{1}\times M_{2}\right)  $ or
$(c,1)(m_{1},m_{2})\in\left(  N_{1}\times M_{2}\right)  $. So, clearly we
conclude $abm_{1}\in N_{1}$ or $cm_{1}\in N_{1}$. Therefore, $N_{1}$ is a
$\psi_{1}-$classical 1-absorbing prime submodule of $M_{1}$. Conversely, let
nonunits $(a_{1},a_{2}),(b_{1},b_{2}),(c_{1},c_{2})\in R$ and $(m_{1}%
,m_{2})\in M$ such that $(a_{1},a_{2})(b_{1},b_{2})(c_{1},c_{2})(m_{1}%
,m_{2})\in N\setminus\phi\left(  N\right)  $. Since $\psi_{2}\left(
M_{2}\right)  =M_{2}$, we get $a_{1}b_{1}c_{1}m_{1}\in N_{1}\setminus\psi
_{1}\left(  N_{1}\right)  $ and this implies that either $abm_{1}\in N_{1}$ or
$cm\in N_{1}$. Thus $(a_{1},a_{2})(b_{1},b_{2})(m_{1},m_{2})\in N$ $\ $or
$(c_{1},c_{2})(m_{1},m_{2})\in N$.
\end{proof}

\begin{theorem}
\label{theo21}Let $R=R_{1}\times R_{2}$ be a decomposable ring and
$M=M_{1}\times M_{2}$ be an $R$-module where $M_{1}$ is an $R_{1}$-module and
$M_{2}$ is an $R_{2}$-module. Suppose $N_{1},N_{2}$ are proper submodules of
$M_{1},M_{2}$, respectively. $\psi_{i}:S(M_{i})\longrightarrow S(M_{i}%
)\cup\left\{  \emptyset\right\}  $ be functions for $i=1,2$ and let $\phi
=\psi_{1}\times\psi_{2}$. If $N=N_{1}\times N_{2}$ is a $\phi-$classical
1-absorbing prime submodule of $M$, then $N_{1}$ is a $\psi_{1}-$classical
1-absorbing prime submodule of $M_{1}$ and $N_{2}$ is a $\psi_{2}-$classical
1-absorbing prime submodule of $M_{2}$.
\end{theorem}

\begin{proof}
Suppose that $N=N_{1}\times N_{2}$ is a $\phi-$ classical 1-absorbing prime
submodule of $M$. Let $abcm\in N_{1}\setminus\psi_{1}\left(  N_{1}\right)  $
that nonunits $a,b,c\in R_{1}$ and $m\in M_{1}$. Get an element $n\in N_{2}$.
We have $(a,1)(b,1)(c,1)(m,n)\in N\setminus\phi(N)$. Then $(a,1)(b,1)(m,n)\in
N$ or $(c,1)(m,n)\in N$. Thus $abm\in N_{1}$ or $cm\in N_{1}$ and thus $N_{1}$
is a $\psi_{1}-$classical 1-absorbing prime submodule of $M_{1}$. By a similar
argument we can show that $N_{2}$ is a $\psi_{2}-$classical 1-absorbing prime
submodule of $M_{2}$.
\end{proof}

\end{document}